\documentclass[12pt, reqno]{amsart}

\usepackage{amssymb}
\usepackage{amsthm}
\usepackage{amsmath}

\newcommand{\F}{{\mathbb F}}

\newcommand{\C}{\mathbb C}

\newcommand{\Z}{\mathbb Z}
\newcommand{\Q}{\mathbb Q}

\newcommand{\modulo}{{\rm mod\ }}
\newcommand{\order}{{\rm ord}}
\newcommand{\trace}{{\rm Tr}}

\newcommand{\cyclic}[1]{\ensuremath{\langle #1 \rangle}}
\newcommand{\zp}[1]{\ensuremath{\Z/#1\Z}}

\newcommand{\gal}{{\rm Gal}}

\newtheorem{theorem}{Theorem}[section]
\newtheorem{lemma}[theorem]{Lemma}
\newtheorem{Rem}[theorem]{Remark}

\newtheorem{coro}[theorem]{Corollary}
\newtheorem{conj}[theorem]{Conjecture}

\numberwithin{equation}{section}

\title[Strongly Regular Cayley Graphs]
{Constructions of Strongly Regular Cayley Graphs using Even Index Gauss Sums}
\author{Fan Wu}
\date{}

\begin{document}
\maketitle

\begin{abstract}
In this paper, generalizing the result in \cite{GXY}, we construct strongly regular Cayley graphs by using union of cyclotomic classes of $\F_q$ and Gauss sums of index $w$, where $w\geq 2$ is even. In particular, we obtain three infinite families of strongly regular graphs with new parameters.
\end{abstract}

\section{Introduction}
We assume that the reader is familiar with the basic theory of strongly regular graphs as can found in \cite{BH, GR, VW}. All graphs considered in this paper are simple and undirected.

A {\it strongly regular graph} srg$(v, k, \lambda, \mu)$ is a regular graph of order $v$ and valency $k$, neither complete nor edgeless, which has the following properties:\\
(1) Any two adjacent vertices have exactly $\lambda$ common neighbors.\\
(2) Any two nonadjacent vertices have exactly $\mu$ common neighbors.\\

Strongly regular graphs have been studied extensively since their introduction by Bose \cite{B} in 1963. In this paper, we will only be concerned with strongly regular Cayley graphs, which are srg with an automorphism group acting sharply transitively on the vertex sets. Such srg are closely related to two-weight linear codes, projective two-intersection sets in finite geometry, and partial difference sets. We refer the reader to \cite{Ma}, \cite{BH}, \cite{SW} for these connections.

Let $\Gamma$ be a graph and $A$ be its adjacency matrix. Then $A$ is a symmetric $(0, 1)$-matrix with all its diagonal entries being zero. The {\it eigenvalues} of $\Gamma$ are by definition the eigenvalues of $A$. Below is a spectral characterization of srg; see \cite[p.~115]{BH} for a proof. For convenience, we call an eigenvalue of $\Gamma$ {\it restricted} if it has an eigenvector orthogonal to the all-one vector. 

\begin{theorem}\label{srgequiv}
For a regular graph $\Gamma$ of order $v$ and valency $k$, not complete nor edgeless, with adjacency matrix $A$, the following are equivalent:\\
(1) $\Gamma$ is an srg$(v, k, \lambda, \mu)$ for certain integers $\lambda, \mu$.\\
(2) $A^{2}=(\lambda-\mu)A+(k-\mu)I+\mu J$, where $I, J$ are the identity matrix and the all-ones matrix, respectively.\\
(3) $A$ has precisely two distinct restricted eigenvalues.
\end{theorem}

%We refer the reader to \cite[p.~115]{BH} for a proof of Theorem \ref{srgequiv}. 

%The two distinct restricted eigenvalues of an srg are usually denoted by $r$ and $s$, where $r$ is the %positive eigenvalue and $s$ the negative one. 

Let $q$ be a prime power, $\F_{q}$ be the finite field of order $q$, and $\F_{q}^{*}=\F_{q}\setminus\{0\}$. Let $D$ be a subset of $\F_{q}^{*}$ such that $-D=D$. The {\it Cayley graph} $\mbox{Cay}(\F_{q}, D)$ is the graph whose vertex set is $\F_{q}$, and two vertices are adjacent if and only if their difference belongs to $D$.  Let $D$ be a subgroup of $\F_{q}^{*}$ such that $-D=D$.  If Cay$(\F_{q}, D)$ is an srg,  then we say that Cay$(\F_{q}, D)$ is a {\it cyclotomic} strongly regular graph. As an example of cyclotomic srg, we mention the Paley graph ${\rm P}(q)$, which is nothing but $\mbox{Cay}(\F_{q}, D)$, where $D$ is the subgroup of $\F_{q}^{*}$ of index 2, and $q$ is a prime power congruent to 1 modulo 4.

Cyclotomic strongly regular graphs have been extensively studied. Let $p$ be a prime, $f$ a positive integer, and $D$ be a subgroup of $\F_{p^{f}}^{*}$ of index $N>1$. If $D$ is the multiplicative group of a subfield of $\F_{p^{f}}$, then it is easy to see that Cay$(\F_{p^{f}}, D)$ is an srg. Such cyclotomic srg are called {\it subfield examples}. Next, if there exists a positive integer $t$ such that $-1\equiv p^{t}\; (\modulo N)$, then it can be shown that Cay$(\F_{q}, D)$ is an srg. (For a proof of this fact, see \cite{BMW}.)  Such cyclotomic srg are called {\it semi-primitive examples}.  Schmidt and White \cite{SW} proposed a conjectural classification of all cyclotomic srg. We state their conjecture below.

\begin{conj} \label{SWConj}
{\em (Conjecture 4.4, \cite{SW})} Let $p$ be a prime, $f$ a positive integer, and $q=p^{f}$. Let $N>1$ be a divisor of $(q-1)/(p-1)$. Assume that $D$ is the subgroup of $\F_{q}^{*}$ of index $N$ such that $-D=D$. If Cay$(\F_{q}, D)$ is an srg, then one of the following holds:
\begin{enumerate}
\item (subfield case) $D=\F_{p^{e}}^{*}$ for some integer $e\geq 1$, $e|f$.
\item (semi-primitive case) There exists a positive integer $t$ such that $-1\equiv p^{t}\;(\modulo N)$.
\item (exceptional case) Cay$(\F_{p^{f}}, D)$ is one of the 11 sporadic examples appearing in the following table:
\begin{eqnarray*}
\begin{tabular}{c c c c}
\hline
$N$ & $p$ & $f$ & $[(\Z/N\Z)^{*}: \cyclic{p}]$\\
\hline
11 & 3 & 5 & 2\\
19 & 5 & 9 & 2\\
35 & 3 & 12 & 2\\
37 & 7 & 9 & 4\\
43 & 11 & 7 & 6\\
67 & 17 & 33 & 2\\
107 & 3 & 53 & 2\\
133 & 5 & 18 & 6\\
163 & 41 & 81 & 2\\
323 & 3 & 144 & 2\\
499 & 5 & 249 & 2\\
\hline
\multicolumn{3}{c}{\quad\quad\quad\quad\quad{\bf Table \quad I}}\\\\
\end{tabular}
\end{eqnarray*}
\end{enumerate}
\end{conj}

Partial results on Conjecture~\ref{SWConj} can be found in \cite{SW, AL}.  However the conjecture as a whole remains open.  On the other hand, in a series recent papers \cite{FX, FMX, GXY}, by using unions of cyclotomic classes of $\F_q$ (instead of using a single cyclotomic class), it is shown that most of the examples in Table I can be generalized into infinite families. Specifically, all index 2 examples but the first one, and the index 4 example in Table I have been generalized into infinite families. The contructions in \cite{FX, FMX, GXY} rely on explicit determination of index 2 or 4 Gauss sums.  A natural question arises: can  constructions similar to those in \cite{FX, FMX, GXY} produce srg in high index cases? In particular, can one generalize the index 6 examples in Table I into infinite families.  A major obstacle is the determination of Gauss sums of high indices. After studying \cite{GXY} closely, we realize that in order to construct strongly regular Cayley graphs by using methods similar to those in \cite{FX, FMX, GXY}, one does not need to evaluate Gauss sums of high indices explicitly;  it suffices to know which subfield (of the cyclotomic field) the Gauss sums belong to.  The results of the paper are as follows. We first generalize the construction of \cite{GXY} to the index $w$ case, where $w\geq 2$ is even. See Theorem~\ref{eigennumber}. This step is straightforward. The main new result is Theorem~\ref{srgcond}, in which we give necessary and sufficient conditions for the construction in Theorem~\ref{eigennumber} to give srg. In Section 4, we use these results to construct explicit families of srg. In particular, we obtain three infinite families of srg. The first family generalizes Example 5 in Table I.
%, and it has parameters
%$$v=11^{7\cdot 43^{m-1}}, \ k=\frac{v-1}{43}, \; r=\frac{21\cdot 11^{\frac{ 7\cdot 43^{m-1}  -1} {2} } -1}
%{43},\; {\rm and}\ s=\frac{-22\cdot 11^{\frac{7\cdot 43^{m-1} -1} {2} } -1}{43},$$
%where $m\geq 1$ is an integer. (Note that the $\lambda$ and $\mu$ values of the srg can be computed %from $v,k,r$ and $s$.) 
The second and third families of srg generalize some subfield examples of cyclotomic strongly regular graphs.

\section{Gauss sums}

Let $p$ be a prime, $f$ be a positive integer, and $q=p^{f}$. Let $\F_{q}$ be the finite field of order $q$. Let $\xi_{p}$ be a fixed complex primitive $p^{\rm th}$ root of unity, and $\trace_{q/p}$ be the trace from $\F_{q}$ to $\F_{p}$. The {\it additive characters} of $\F_{q}$ are the homomorphisms from the additive group $(\F_{q}, +)$ to $\C^{*}$, the multiplicative group of $\C$, and they are given by
\begin{equation*}
\psi_{a}(x)=\xi_{p}^{\trace_{q/p}(ax)},
\end{equation*}
where $a\in \F_{q}$. We usually write $\psi_{1}$ simply as $\psi$, which is called the {\it canonical} additive character of $\F_{q}$. The {\it multiplicative characters} of $\F_{q}$ are the homomorphisms from the multiplicative group $\F_{q}^{*}$ to $\C^{*}$.

Let $N$ be a positive integer with $N\mid (q-1)$, and let $\chi$ be a multiplicative character of $\F_{q}$ of order $N$. The {\it Gauss sum} $g(\chi)$ of order $N$ is defined by
\begin{equation*}
g(\chi)=\sum_{a\in \F_{q}^{*}}\chi(a)\psi(a).
\end{equation*}
Clearly Gauss sums of order $N$ belong to  $\Z[\xi_{N}, \xi_{p}]$, the integer ring of $\Q(\xi_{N}, \xi_{p})$, where $\xi_N$ is a primitive $N^{\rm th}$ root of unity. Let $\sigma_{a, b}$ be the Galois automorphism of $\Q(\xi_{N}, \xi_{p})$ defined by
\begin{equation*}
\sigma_{a, b}(\xi_{N})=\xi_{N}^{a},\,\, \sigma_{a, b}(\xi_{p})=\xi_{p}^{b}
\end{equation*}
The following lemma gives some useful properties of Gauss sums. For a proof of the lemma, we refer the reader to \cite[p.~10]{BEW} and \cite[p.~208]{IR}.

\begin{lemma}\label{GSumProps}
Let $\chi$ be a multiplicative character of $\F_{q}$ of order $N$. Then
\begin{enumerate}
\item $g(\chi)=-1$, if $\chi=\chi_{0}$ (the trivial character), and $|g(\chi)|=q$, if $\chi\neq \chi_{0}$.
\item $\sigma_{a, b}(g(\chi))=\overline{\chi}^{a}(b)g(\chi^{a})$, where $\overline{\chi}=\chi^{-1}$.
\item $\overline{g(\chi)}=\chi(-1)g(\overline{\chi})$, and $\sigma_{p, 1}(g(\chi))=g(\chi^{p})=g(\chi)$.
\item For a character of $\chi$ of order $N$, $g(\chi)^{N}\in \Z[\xi_{N}]$.
\end{enumerate}
\end{lemma}

For the rest of this paper, we assume that (i) $\gcd(p(p-1), N)=1$, where $N|(q-1)$, and $q=p^f$, $f$ is the order of $p$ modulo $N$, (ii) $-1\not\in \langle p\rangle$, the cyclic subgroup of  $(\Z/N\Z)^{*}$ generated by $p$. These assumptions have the following consequences. 

First, the index $[(\Z/N\Z)^{*}: \cyclic{p}]$, denoted by $w$, must be even. This can be seen as follows. From $\gcd(p(p-1), N)=1$, we see that $N$ is odd. Thus $\phi(N)$ is even, where $\phi$ is the Euler totient function. If $w$ is odd, then $f=\phi(N)/w$ is even. It follows that $p^{f/2}\equiv -1 \pmod N$, contradicting the assumption that $-1\not\in \cyclic{p}$.

Secondly, let $\chi$ be a multiplicative character of $\F_q$ of order $N$. We claim that $g(\chi)\in \Z[\xi_{N}]$. We prove this claim as follows. For any $b\in \F_p^*$, since $\gcd(N, p-1)=1$, we have $\chi(b)=1$. Hence by Part (2) of Lemma~\ref{GSumProps},  $\sigma_{1,b}(g(\chi))={\bar \chi}(b)g(\chi)=g(\chi)$. It follows that $g(\chi)\in \Z[\xi_{N}]$. We can actually go a little further. Let $K$ be the decomposition field of the prime $p$ in $\Q(\xi_{N})$. Then it is well known \cite[p.~197]{IR} that $\gal(\Q(\xi_{N})/K)=\langle \sigma_{p,1}\rangle$. By Part (3) of Lemma~\ref{GSumProps}, we have  $g(\chi)\in K$. In fact, we have $g(\chi)\in O_K$, the integer ring of $K$.

Gauss sums $g(\chi)$ with $\chi$ being a multiplicative character of order $N$ of $\F_q$ (and $-1\not\in\langle p\rangle$) are called {\it Gauss sums of index $w$}.

For the remainder of this paper we will further assume that $N=p_1^m$, where $p_1$ is an odd prime, $m\geq 1$ is an integer, and $w|(p_1 -1)$. In this case, $\gal(\Q(\xi_{N})/\Q)\cong (\Z/N\Z)^{*}$ is cyclic, and $K$ is the unique imaginary subfield of $\Q(\xi_{N})$ with $[K:\Q]=w$. Since $w|(p_1 -1)$, we see that $K$ is in fact a subfield of $\Q(\xi_{p_{1}})$. Therefore if $\chi$ is a multiplicative character of $\F_q$ of order $N$, we in fact have $g(\chi)\in \Z[\xi_{p_{1}}]$.

Let $g$ be a primitive root modulo $p_1$. Define $\widetilde{C_{j}}=g^{j}\cyclic{p}\subseteq (\zp{p_{1}})^{*}$, for all $0\le j\le w-1$. Then $ (\zp{p_{1}})^{*}=\cup_{j=0}^{w-1} \widetilde{C_{j}}.$ For $0\le j\le w-1$, we define $\eta_{j}$ by
\begin{equation}\label{Gaussperiods}
\eta_{j}=\sum_{a\in \tilde{C_{j}}}\xi_{p_{1}}^{a},
\end{equation}
where $\xi_{p_{1}}$ is a complex primitive $p_{1}$-th root of unity. The following lemma is well known (see \cite{T}).

\begin{lemma}\label{IntegralBasis}
With the above assumptions, $\{\eta_{j}\mid 0\le j\le w-1\}$ is an integral basis of $K$.
\end{lemma}

Let $\chi$ be a multiplicative character of $\F_q$ of order $N$. Let $r$ be the largest nonnegative integer such that $p^r|g(\chi)$. That is, $p^{-r}g(\chi)\in O_K$, but $p^{-(r+1)}g(\chi)\not\in O_K$. Note that if $\chi'$ is another multiplicative character of $\F_q$ of order $N$, then there exists a $d\in (\Z/N\Z)^*$ such that $\chi'=\chi^d$; it follows that $g(\chi')=g(\chi^d)=\sigma_{d,1}(g(\chi))$. This shows that the integer $r$ does not depend on the choice of the multiplicative character of order $N$. The explicit computation of $r$ can be done by using Stickelberger's theorem on the prime ideal factorization of Gauss sums. We state the following lemma,  whose proof can be found in \cite{FYL}.

\begin{lemma}\label{pPower}
Let $\chi$ be a a multiplicative character of $\F_q$ of order $N$. With the above assumptions and notation, we have
$$r=\frac{f -\tilde{f}}{2} +b,$$
where $\tilde{f}=\frac{p_1 -1}{w}$, $b={\rm min}\{b_0, b_1, \ldots ,b_{w-1}\}$ and $b_j=\frac{1}{p_1} \sum _{z \in ([1,p_1-1]\cap \tilde{C}_j)} z$ for all $0\leq j \leq w-1$, here $[1,p_1-1]$ denotes the set of integers $x$, $1\leq x\leq p_1-1$.
\end{lemma}

\section{Cyclotomic classes and strongly regular Cayley graphs}

In this section, $N$, $p_{1}$, $p$, $f$ and $w$ are the same as those in Section 2. Let $\F_{q}$ be the finite field of order $q$, $q=p^{f}$. Let $\gamma$ be a fixed primitive element of $\F_{q}$. The $N^{\rm th}$  {\it cyclotomic classes} $C_{0}, C_{1}, \ldots, C_{N-1}$ of $\F_q$ are defined by
\begin{equation*}
C_{i}:=\{\gamma^{i+jN}\mid 0\le j\le \frac{q-1}{N}-1\},
\end{equation*}
where $0\le i\le N-1$. Clearly $C_{i}=\gamma^{i}C_{0}$, for $0\leq i\leq N-1$. Let $\psi$ be the canonical additive character of $\F_{q}$. The $N^{\rm th}$ {\it cyclotomic periods} (also known as {\it Gauss periods}) are defined by
\begin{equation}\label{Gaussperiodff}
\tau_{a}:=\sum_{x\in C_{a}}\psi(x),
\end{equation}
where $0\le a\le N-1$.

Let $\widehat{\F_{q}^{*}}$ be the group of multiplicative characters of $\F_{q}$. Gauss sums can be viewed as the Fourier coefficients in the Fourier expansion of the additive characters in terms of the multiplicative characters of $\F_{q}$. That is,
\begin{equation}\label{AddChar}
\psi(a)=\frac{1}{q-1}\sum_{\chi\in \widehat{\F_{q}^{*}}}g(\overline{\chi})\chi(a), \mbox{ for all } a\in \F_{q}^{*}
\end{equation}
where $\overline{\chi}=\chi^{-1}$. Using (\ref{AddChar}), we obtain the relationship between Gauss sums and Gauss periods.

\begin{lemma}\label{Gaussperiodcomp}
For $0\leq a\leq N-1$, we have 
\begin{equation}\label{}
\tau_{a}=\frac{1}{N}\sum_{\chi\in C_{0}^{\perp}}g(\overline{\chi})\chi(\gamma^{a}),
\end{equation}
where $C_0^{\perp}$ is the subgroup of $\widehat{\F_{q}^{*}}$ consisting of characters which are trivial on $C_0$.
\end{lemma}
\noindent{\bf Proof}: From (\ref{Gaussperiodff}) and (\ref{AddChar}), we have
\begin{eqnarray*}
\tau_{a} &=& \sum_{x\in C_{0}}\psi(\gamma^{a}x)=\frac{1}{q-1}\sum_{x\in C_{0}}\sum_{\chi\in \widehat{\F_{q}^{*}}}g(\overline{\chi})\chi(\gamma^{a}x)\\
&=& \frac{1}{q-1}\sum_{\chi\in \widehat{\F_{q}^{*}}}\sum_{x\in C_{0}}g(\overline{\chi})\chi(\gamma^{a}x)\\
&=&  \frac{1}{q-1}\sum_{\chi\in \widehat{\F_{q}^{*}}}g(\overline{\chi})\chi(\gamma^{a})\sum_{x\in C_{0}}\chi(x)
\end{eqnarray*}

If $\chi\not\in C_{0}^{\perp}$, then there exists a positive integer $\ell$ such that $\chi(\gamma^{\ell N})\neq 1$. We have $$\chi(\gamma^{\ell N})\sum_{x\in C_{0}}\chi(x)=\sum_{x\in C_{0}}\chi(\gamma^{\ell N}x)=\sum_{x\in C_{0}}\chi(x).$$ It follows that $\sum_{x\in C_{0}}\chi(x)=0$ since $\chi(\gamma^{\ell N})\neq 1$. If $\chi\in C_{0}^{\perp}$, clearly $\sum_{x\in C_{0}}\chi(x)=|C_{0}|=(q-1)/N$.  The proof of the lemma is now complete. $\hfill\Box$

Define
\begin{equation}\label{DefD}
D:=\bigcup_{i=0}^{p_{1}^{m-1}-1}C_{i}.
\end{equation}
Using $D$ as connection set, we construct ${\rm Cay}(\F_q, D)$.

\begin{theorem}\label{eigennumber}
The Cayley graph Cay($\F_{q}$, $D$) is an undirected regular graph of valency $|D|$. It has at most $w+1$ distinct restricted eigenvalues.
\end{theorem}

The proof is parallel to that of Theorem 3.1 in \cite{GXY}. Since we will use some parts of the proof later on, we will give the complete proof here.

\noindent{\bf Proof}: Since $N=p_1^m$ is odd, we have $2N |(q-1)$ or $q$ is even; consequently $-C_0=C_0$. Hence $-D=D$. It follows that Cay($\F_{q}$, $D$) is undirected. Now $0\not\in D$, we see that Cay($\F_{q}$, $D$) has no loops. From definition, we deduce that Cay($\F_{q}$, $D$) is a regular graph of valency $|D|$. 

The restricted eigenvalues of Cay($\F_{q}$, $D$), as explained in \cite[p.~136]{BH}, are given by 
\begin{equation*}
\psi(\gamma^{a}D)=\sum_{x\in D}\psi(\gamma^{a}x),
\end{equation*}
where $0\le a\le N-1$. By (\ref{DefD}) and Lemma \ref{Gaussperiodcomp}, we have
\begin{eqnarray}\label{eigenguass}
\psi(\gamma^{a}D) &=& \sum_{i=0}^{p_{1}^{m-1}-1}\psi(\gamma^{a}C_{i}) = \sum_{i=0}^{p_{1}^{m-1}-1}\tau_{i+a}\nonumber\\
&=& \frac{1}{N}\sum_{\chi\in C_{0}^{\perp}}g(\overline{\chi})\sum_{i=0}^{p_{1}^{m-1}-1}\chi(\gamma^{a+i}).
\end{eqnarray}
If $\chi\in C_0^{\perp}$ and $\chi=\chi_0$ (the trivial character), then $g(\overline{\chi})=-1$ and $\sum_{i=0}^{p_{1}^{m-1}-1}\chi(\gamma^{a+i})=p_{1}^{m-1}$. If $\chi\in C_0^{\perp}$ and $\order(\chi)\neq 1$, then $\order(\chi)=p_{1}^{\ell}$, $1\le \ell\le m$ since $\order(\chi)||C_{0}^{\perp}|$. For those characters $\chi$ with $1\neq \order(\chi)<p_1^{m}$, we have
\begin{eqnarray*}
\sum_{i=0}^{p_{1}^{m-1}-1}\chi(\gamma^{a+i})=\chi(\gamma^{a})\sum_{i=0}^{p_{1}^{m-1}-1}\chi(\gamma)^{i}=\chi(\gamma^{a})\frac{\chi(\gamma)^{p_{1}^{m-1}}-1}{\chi(\gamma)-1}=0.
\end{eqnarray*}
Thus, in (\ref{eigenguass}), the terms corresponding to characters of order $p_{1}^{\ell}$, $1\le \ell\le m-1$, vanish. Hence (\ref{eigenguass}) can be simplified to
\begin{equation}\label{eigensimplify}
\psi(\gamma^{a}D)=\frac{1}{N}[-p_{1}^{m-1}+\sum_{\chi\in C_{0}^{\perp}, \order(\chi)=p_{1}^{m}}g(\overline{\chi})\sum_{i=0}^{p_{1}^{m-1}-1}\chi(\gamma^{a+i})].
\end{equation}
Define a multiplicative character $\theta$ of $\F_{q}$ by setting $\theta(\gamma)=\xi_{N}$. Then $\cyclic{\theta}=C_{0}^{\perp}$ since $C_{0}^{\perp}$ is the unique subgroup of $\widehat{\F_{q}^{*}}$ of order $N$. Thus any multiplicative character $\chi$ of order $p_{1}^{m}$ can be expressed as $\theta^{d}$ for some $d$ in $(\Z/N\Z)^{*}$. We have
\begin{eqnarray}\label{eigenprim}
\psi(\gamma^{a}D)=\frac{1}{N}[-p_{1}^{m-1}+\sum_{d\in (\Z/N\Z)^{*}}g(\overline{\theta}^{d})\sum_{i=0}^{p_{1}^{m-1}-1}\theta^{d}(\gamma^{a+i})]
\end{eqnarray}
For convenience, we set
\begin{equation}\label{sa}
S_{a}:=\sum_{d\in (\Z/N\Z)^{*}}g(\overline{\theta}^{d})\sum_{i=0}^{p_{1}^{m-1}-1}\theta^{d}(\gamma^{a+i}).
\end{equation}

Let $r$ be the positive integer given in Lemma~\ref{pPower} such that $p^{-r}g(\overline{\theta})\in O_{K}$. By Lemma \ref{IntegralBasis}, we have 
\begin{equation}\label{expressionGsum}
g(\overline{\theta})=p^{r}(N_{0}\eta_{0}+\cdots + N_{w-1}\eta_{w-1}),
\end{equation} 
where $N_{0}, \ldots, N_{w-1}$ are integers and $\eta_{0}, \ldots, \eta_{w-1}$ are defined in (\ref{Gaussperiods}). From Lemma \ref{GSumProps}, we have $g(\overline{\theta}^{d})=\sigma_{d,\;1}(g(\overline{\theta}))$. To simply notation, we simply write $\sigma_d$ for $\sigma_{d,1}$. It follows that
\begin{eqnarray*}
g(\overline{\theta}^{d})&=&\sigma_{d}(g(\overline{\theta}))=\sigma_{d}(p^{r}(N_{0}\eta_{0}+\cdots + N_{w-1}\eta_{w-1}))\nonumber\\
&=& p^{r}(N_{0}\eta_{0}^{\sigma_{d}}+\cdots + N_{w-1}\eta_{w-1}^{\sigma_{d}}).
\end{eqnarray*}
Now writing $d\in (\Z/N\Z)^{*}$ as $d=d_{1}+p_{1}d_{2}$, where $d_{1}\in (\Z/p_{1}\Z)^{*}$ and $d_{2}\in \Z/p_{1}^{m-1}\Z$, we have $\eta_{j}^{\sigma_{d}}=\sigma_{d}(\sum_{c\in \tilde{C_{j}}}\xi_{p_{1}}^{c})=\sum_{c\in \tilde{C_{j}}}\sigma_{d_{1}+p_{1}^{m-1}d_{2}}(\xi_{p_{1}}^{c})=\eta_{j}^{\sigma_{d_{1}}}$. We have
\begin{eqnarray}\label{uniexp}
g(\overline{\theta}^{d})&=& p^{r}(N_{0}\eta_{0}^{\sigma_{d}}+\cdots+N_{w}\eta_{w}^{\sigma_{d}})\nonumber\\
&=& p^{r}(N_{0}\eta_{0}^{\sigma_{d_{1}}}+\cdots+N_{w}\eta_{w}^{\sigma_{d_{1}}})\nonumber\\
&=& \sigma_{d_{1}}(p^{r}(N_{0}\eta_{0}+\cdots+N_{w}\eta_{w}))\nonumber\\
&=& \sigma_{d_{1}}(g(\overline{\theta}))=g(\overline{\theta}^{d_{1}}).
\end{eqnarray}
Hence (\ref{sa}) can be written as
\begin{eqnarray*}
S_{a} &=& \sum_{d\in (\Z/N\Z)^{*}}g(\overline{\theta}^{d})\sum_{i=0}^{p_{1}^{m-1}-1}\theta^{d}(\gamma^{a+i})\\
&=& \sum_{d_{1}\in (\Z/p_{1}\Z)^{*}}\sum_{d_{1}\in \Z/p_{1}^{m-1}\Z}g(\overline{\theta}^{d_{1}+p_{1}d_{2}})\sum_{i=0}^{p_{1}^{m-1}-1}\theta^{d_{1}+p_{1}d_{2}}(\gamma^{a+i})\\
&=& \sum_{d_{1}\in (\Z/p_{1}\Z)^{*}}\sum_{i=0}^{p_{1}^{m-1}-1}g(\overline{\theta}^{d_{1}})\theta^{d_{1}}(\gamma^{a+i})\sum_{d_{2}\in \Z/p_{1}^{m-1}\Z}\theta^{p_{1}d_{2}}(\gamma^{a+i}).
\end{eqnarray*}
Note that $\sum_{d_{2}\in \Z/p_{1}^{m-1}\Z}\theta^{d_{2}p_{1}(a+i)}(\gamma)=0$ if and only if $p_{1}^{m-1}\nmid (a+i)$. We only need to consider the terms for which $p_{1}^{m-1}\mid (a+i)$. For each $0\leq a\leq N-1$, there exists a unique $i_a\in \{0,1,\ldots ,p_1^{m-1}-1\}$ such that $p_{1}^{m-1}\mid (a+i_a)$; write $a+i_a=p_{1}^{m-1}j_{a}$, $j_{a}\in \Z/p_{1}\Z$, we have
\begin{equation}\label{sa1}
S_{a}=p_{1}^{m-1}\sum_{d_{1}\in (\Z/p_{1}\Z)^{*}}g(\overline{\theta}^{d_{1}})\theta^{d_{1}}(\gamma^{p_{1}^{m-1}j_{a}})=p_{1}^{m-1}\sum_{d_{1}\in (\Z/p_{1}\Z)^{*}}g(\overline{\theta}^{d_{1}})\xi_{p_{1}}^{j_{a}d_{1}}.
\end{equation}
For each $j\in \Z/p_1\Z$, define an additive character $\psi_{j}$ on $\Z/p_{1}\Z$ such that $\psi_{j}(d_{1})=\xi_{p_{1}}^{jd_{1}}$. We have
\begin{eqnarray}\label{sa3}
S_{a} &=& p_{1}^{m-1}\sum_{d_{1}\in (\Z/p_{1}\Z)^{*}}g(\overline{\theta}^{d_{1}})\psi_{j_{a}}(d_{1})\nonumber\\
&=& p_{1}^{m-1}p^{r}\sum_{i=0}^{w-1}\sum_{d_{1}\in \widetilde{C_{i}}}(N_{0}\eta_{0}^{\sigma_{d_{1}}}+\cdots+N_{w-1}\eta_{w-1}^{\sigma_{d_{1}}})\psi_{j_{a}}(d_{1})\nonumber\\
&=& p_{1}^{m-1}p^{r}[(N_{0}\eta_{0}+N_{1}\eta_{1}+\cdots+N_{w-1}\eta_{w-1})\sum_{d_{1}\in \widetilde{C_{0}}}\psi_{j_{a}}(d_{1})\nonumber\\
&&+(N_{0}\eta_{1}+N_{1}\eta_{2}+\cdots+N_{w}\eta_{0})\sum_{d_{1}\in \widetilde{C_{1}}}\psi_{j_{a}}(d_{1})\nonumber\\
&& \cdots\nonumber\\
&& +(N_{0}\eta_{w-1}+N_{1}\eta_{0}+\cdots+N_{w-1}\eta_{w-2})\sum_{d_{1}\in \widetilde{C_{w-1}}}\psi_{j_{a}}(d_{1})].
\end{eqnarray}

Let $M_{0}=N_{0}+N_{1}+\cdots+N_{w-1}$. Note that $\sum_{i=0}^{w-1}\eta_{i}=-1$. We continue the computations of $S_a$ by considering two case.

\noindent {\bf Case 1.} $j_{a}=0$.  In this case, $\psi_{j_{a}}(d_{1})=1$ for all $d_1\in \Z/p_1\Z$. It follows that $\sum_{d_{1}\in \widetilde{C_{z}}}\psi_{j_{a}}(d_{1})=(p_{1}-1)/w$, $0\le z\le w-1$. Thus we have 
\begin{equation*}\label{sajais0}
S_{a}=\frac{p_{1}-1}{w}p_{1}^{m-1}p^{r}(N_{0}\sum_{i=0}^{w-1}\eta_{i}+N_{1}\sum_{i=0}^{w-1}\eta_{i}+\cdots+N_{w-1}\sum_{i=0}^{w-1}\eta_{i})
=\frac{1-p_{1}}{w}p_{1}^{m-1}p^{r}M_{0}
\end{equation*}

\noindent {\bf Case 2.} $j_{a}\neq 0$. In this case, $j_{a}$ must belong to a unique coset of $\cyclic{p}$ in $(\Z/p_{1}\Z)^{*}$, say $j_{a}\in g^{t}\cyclic{p}$, where $0\le t\le w-1$. In this case,  for any $0\le z\le w-1$, we have
\begin{equation*}
\sum_{d_{1}\in \widetilde{C_{z}}}\psi_{j_{a}}(d_{1})=\eta_{\overline{z+t}},
\end{equation*}
where $\overline{z+t}$ is the least nonnegative residue of $z+t$ modulo $w$.

Define
\begin{eqnarray*}\label{corrsum}
\begin{aligned}
& K_{0} = \eta_{0}^{2}+\cdots+\eta_{w-1}^{2},\\
& K_{1} = \eta_{0}\eta_{1}+\cdots+\eta_{w-1}\eta_{0},\\
& \cdots\\
& K_{w-1} = \eta_{0}\eta_{w-1}+\cdots+\eta_{w-1}\eta_{w-2}.
\end{aligned}
\end{eqnarray*}
Then
\begin{eqnarray}\label{sajanot0}
S_{a} &=& p_{1}^{m-1}p^{r}[(N_{0}\eta_{0}+N_{1}\eta_{1}+\cdots+N_{w-1}\eta_{w-1})\eta_{\overline{t}}\nonumber\\
&&+(N_{0}\eta_{1}+N_{1}\eta_{2}+\cdots+N_{w-1}\eta_{0})\eta_{\overline{1+t}}\nonumber\\
&& \cdots\nonumber\\
&& +(N_{0}\eta_{w-1}+N_{1}\eta_{0}+\cdots+N_{w-1}\eta_{w-2})\eta_{\overline{w-1+t}}]\nonumber\\
&=& p_{1}^{m-1}p^{r}[N_{0}(\eta_{0}\eta_{\overline{t}}+\eta_{1}\eta_{\overline{1+t}}+\cdots+\eta_{w-1}\eta_{\overline{w-1+t}})\nonumber\\
&&+N_{1}(\eta_{0}\eta_{\overline{w-1+t}}+\eta_{1}\eta_{\overline{t}}+\cdots+\eta_{w-1}\eta_{\overline{w-2+t}})\nonumber\\
&& \cdots\nonumber\\
&& +N_{w-1}(\eta_{0}\eta_{\overline{1+t}}+\eta_{1}\eta_{\overline{2+t}}+\cdots+\eta_{w-1}\eta_{\overline{t}})]\nonumber\\
&=& p_{1}^{m-1}p^{r}(N_{0}K_{\overline{t}}+N_{1}K_{\overline{t+1}}+\cdots+N_{w-1}K_{\overline{t+w-1}})
\end{eqnarray}

In Section 2, we have shown that $w$ is even and $f$ is odd. The Gauss periods $\eta_j$ satisfy the following relations (see \cite{T}):
$$K_{w/2}=(1+(w-1)p_{1})/w,  \; K_{j}=(1-p_{1})/w, \;{\rm if}\; j\neq w/2.$$

Clearly there exists a unique element in the set $\{\overline{t},\ \overline{t+1}, \ldots, \overline{t+w-1}\}$ that is equal to $w/2$, say $\overline{t+h(a)}=w/2$, where $0\leq h(a)\leq w-1$, and $h(a)$ depends on $a$. We have
\begin{eqnarray}\label{sajanot0final}
S_{a} &=& p_{1}^{m-1}p^{r}[\frac{1-p_{1}}{w}(M_{0}-N_{h(a)})+\frac{1+(w-1)p_{1}}{w}N_{h(a)}]\nonumber\\
&=& p_{1}^{m-1}p^{r}(\frac{1-p_{1}}{w}M_{0}+p_{1}N_{h(a)}).
\end{eqnarray}
Therefore in this case, we have
\begin{equation*}
S_{a}\in \{p_{1}^{m-1}p^{r}(\frac{1-p_{1}}{w}M_{0}+p_{1}N_{i})\mid 0\le i\le w-1\}.
\end{equation*}

Summing up, let
\begin{eqnarray}\label{eignset}
E&=&\left\{\frac{1}{p_{1}}(-1+\frac{1-p_{1}}{w}p^{r}M_{0})\right\}\nonumber\\
&&\cup \left\{\frac{1}{p_{1}}(-1+\frac{1-p_{1}}{w}p^{r}M_{0})+p^{r}N_{i}\mid 0\le i\le w-1\right\}.
\end{eqnarray}
We have shown that the restricted eigenvalues of Cay$(\F_q,D)$ belong to $E$. Since $|E|\leq w+1$, we see that Cay$(\F_q,D)$ has at most $w+1$ distinct restricted eignevalues. The proof of the theorem is now complete. $\hfill\Box$

Next we give necessary and sufficient conditions for ${\rm Cay}(\F_q, D)$  to be an srg. The proof uses discrete Fourier transforms, which were first employed in the proof of Theorem~3.1 in \cite{SW}.

\begin{theorem}\label{srgcond}
Let $p_1$ be a prime, $m\geq 1$, $N=p_1^m$. Let $p\neq p_1$ be a prime, $f={\rm ord}_N(p)$, $w=\phi(N)/f$, and $q=p^f$. Assume that $-1\not\in\langle p\rangle$, $\gcd(p(p-1), N)=1$ and $w|(p_1-1)$. Let $r$ be given in Lemma~\ref{pPower} and $D$ be defined as in (\ref{DefD}). Then Cay($\F_{q}$, $D$) is a strongly regular graph if and only if there exists an integer $\ell$, $1\leq \ell\leq w-1$, such that
\begin{equation}\label{cond1}
\frac{p^r(1-p_1)\ell}{w}\equiv \epsilon\; (\modulo p_{1})
\end{equation}
and
\begin{equation}\label{cond2}
p^{s}=\frac{\ell}{w}\left(p_{1}-\frac{(p_{1}-1)\ell}{w}\right),
\end{equation}
where $s=f-2r$ and $\epsilon=\pm 1$.
\end{theorem}

\noindent{\bf Proof}: Suppose that Cay($\F_{q}$, $D$) is a strongly regular graph. Then by Theorem~\ref{srgequiv}, Cay($\F_{q}$, $D$) has exactly two distinct restricted eigenvalues. By our computations of the restricted eigenvalues of Cay($\F_{q}$, $D$) in the proof of Theorem~\ref{eigennumber}, we must have $N_{i}\in \{0, \epsilon\}$ for all $0\le i\le w-1$, where $\epsilon\neq 0$ is an integer. Thus (\ref{expressionGsum}) becomes
\begin{equation*}
g(\overline{\theta})=\epsilon p^{r}\sum_{i\in I}\eta_{i},
\end{equation*}
where $I=\{i\mid N_i=\epsilon, 0\leq i\leq w-1\}$. From $|g(\overline{\theta})|^2=p^{f}$, we obtain that 
\begin{equation}\label{abs}
|\sum_{i\in I}\eta_{i}|^{2}=p^{f-2r}/\epsilon^{2}.
\end{equation}
 It follows that $\epsilon$ must be a power of $p$. Since $r$ is the largest power of $p$ dividing the Gauss sum $g(\overline{\theta})$ (see Lemma~\ref{pPower}), we have $\epsilon=\pm 1$.

Let $s=f-2r$ and $D'=\cup_{i\in I}\widetilde{C_{i}}\subset (\Z/p_{1}\Z)^{*}$. From (\ref{abs}), we see that $D'$ is a difference set in $(\Z/p_1\Z, +)$ with parameters $(p_1, \frac{p_1-1}{w}\ell, \frac{p_1-1}{w}\ell-p^s)$, 
where $\ell=|I|$. From the basic parameter relation of difference sets, we obtain that 
$$p^{s}=(\ell/w)(p_{1}-(p_{1}-1)\ell/w).$$

Next we claim that $\frac{p^r(1-p_1)\ell}{w}\equiv \epsilon$ (mod $p_1$). This can be seen as follows.
\begin{eqnarray*}
\psi(\gamma^{a}D) &=& \frac{1}{p_{1}^{m}}(-p_{1}^{m-1}+S_{a})\\
&=& \frac{1}{p_{1}}(-1+\epsilon p^{r}(1-p_{1})\ell/w)+p^{r}N_{h(a)}.
\end{eqnarray*}
Since $\psi(\gamma^{a}D)$ are integers for all $0\leq a\leq N-1$, we see that $(1-p_{1})\ell p^{r}/w\equiv \epsilon \pmod {p_{1}}$. 

Conversely, let
\begin{equation*}
x=\frac{1}{p_{1}}\left(-1+\frac{1-p_{1}}{w} p^{r}\ell \epsilon\right),
\end{equation*}
with $\epsilon=\pm 1$. By (\ref{cond1}), we have $x\in \Z$. Define a function $\varphi: (\Z/N\Z, +)\rightarrow \Z$ by
\begin{equation}\label{Func}
\varphi(a):=\frac{\psi(\gamma^{a}D)-x}{p^{r}}, \; \forall a\in \Z/N\Z.
\end{equation}
We note that since $\psi(\gamma^a D)$ are algebraic integers, and by the computations in the proof of Theorem~\ref{eigennumber}, $\psi(\gamma^a D)=\frac{1}{p_{1}}(-1+\frac{1-p_{1}}{w}p^{r}M_{0})$ or $\frac{1}{p_{1}}(-1+\frac{1-p_{1}}{w}p^{r}M_{0})+p^rN_{h(a)}$, which are rationals, we must have $\psi(\gamma^a D)\in \Z$. It follows that $\frac{1-p_{1}}{w}p^{r}M_{0}\equiv 1\pmod {p_1}$. Now by assumption, we have $\frac{(1-p_1)p^r\ell}{w}\equiv \epsilon\; (\modulo p_{1})$. Therefore $$M_0\equiv \epsilon \ell \pmod{p_1}.$$
We thus have $\varphi(a)\in \Z$ indeed. 

To simplify notation, we use $G$ to denote the cyclic group $(\Z/N\Z, +)$. Then $\widehat{G}=\{\nu^j\mid 0\leq j\leq N-1\}$, where $\nu$ is the character of $G$ sending $1$ to $\xi_N$. The Fourier transform $\hat{\varphi}$ of $\varphi$ is given by
\begin{equation*}
\hat{\varphi}(\nu^{j})=\frac{\sum_{a\in G}\varphi(a)\nu^j(a)}{\sqrt{N}},
\end{equation*}
for $0\le j\le N-1$.

When $j=0$, we have
\begin{eqnarray}\label{FTtrivial}
\hat{\varphi}(\nu^0) &=& \frac{\sum_{a\in G}(\psi(\gamma^{a}D)-x)}{\sqrt{N}p^{r}}=\frac{\sqrt{N}(-1-p_1x)}{p_1p^r}=\sqrt{N}\frac{\epsilon (p_{1}-1)\ell}{wp_{1}}
\end{eqnarray}
For $1\leq j\leq N-1$, we have
\begin{eqnarray}\label{FTnontrivial}
\hat{\varphi}(\nu^{j}) &=&\frac{\sum_{a\in G}(\psi(\gamma^{a}D)-x)\nu^{j}(a)}{\sqrt{N}p^{r}}\nonumber\\
&=&\frac{\sum_{a\in G}\psi(\gamma^{a}D)\nu^{j}(a)}{\sqrt{N}p^{r}}\nonumber\\
\end{eqnarray}
By (\ref{eigenprim}), we have
\begin{eqnarray}\label{final}
\hat{\varphi}(\nu^{j})&=& \frac{1}{p^{r}\sqrt{N}}\sum_{a\in G}\frac{1}{N}\left(-p_{1}^{m-1}+S_a\right)\nu^{j}(a)\nonumber\\
&=& \frac{1}{p^{r}N\sqrt{N}}\sum_{a\in G}S_a\nu^{j}(a)\nonumber\\
&=& \frac{1}{p^{r}N\sqrt{N}} \sum_{d\in (\Z/N\Z)^{*}}g(\overline{\theta}^{d})\sum_{i=0}^{p_1^{m-1}-1}\theta(\gamma)^{di}\sum_{a\in G}\xi_{N}^{a(d+j)}
\end{eqnarray}

If $p_1|j$, then the (inner) sum $\sum_{a\in G}\xi_{N}^{a(d+j)}=0$ since $d\in (\Z/N\Z)^{*}$ (i.e., $d$ is relatively prime to $N$); we thus have $\hat{\varphi}(\nu^{j})=0$.

If $\gcd(p_1,j)=1$, then the (inner) sum $\sum_{a\in G}\xi_{N}^{a(d+j)}$ is nonzero (and equals $N$) if and only if $j\equiv -d \pmod N$; in this case, we have
$$\hat{\varphi}(\nu^{j})=\frac{1}{p^r\sqrt{N}}g(\theta^j)(\sum_{i=0}^{p_1^{m-1}-1}\xi_N^{-ji}).$$

Note that the above formula also holds true for those $j$ such that $1\leq j\leq N-1$ and $p_1|j$ since $\sum_{i=0}^{p_1^{m-1}-1}\xi_N^{-ji}=0$ if $p_1|j$. Therefore for all $1\leq j\leq N-1$, we have
\begin{eqnarray}\label{FTnontrivial1}
\hat{\varphi}(\nu^{j})=\frac{1}{p^r\sqrt{N}}g(\theta^j)(\sum_{i=0}^{p_1^{m-1}-1}\xi_N^{-ji}).
\end{eqnarray}

Using the definition of $\varphi$, we have
\begin{eqnarray}\label{sumfunc1}
\sum_{a\in G}\varphi(a)=\sum_{a\in G}\frac{\psi(\gamma^{a}D)-x}{p^{r}}=N\frac{-1-xp_{1}}{p^{r}p_1}=\frac{N}{p_1}\cdot \frac{\epsilon (p_{1}-1)\ell}{w}.
\end{eqnarray}
From (\ref{FTtrivial}), (\ref{FTnontrivial1}) and Parseval's identity, we have
\begin{eqnarray}\label{sumfunc2}
\sum_{a\in G}\varphi(a)^{2}&=& \sum_{j=0}^{N-1}|\hat{\varphi}(\nu^j)|^{2}=|\hat{\varphi}(\nu^0)|^{2}+\sum_{j=1}^{N-1}|\hat{\varphi}(\nu^j)|^{2}\nonumber\\
&=&\frac{N(p_{1}-1)^2\ell^2}{w^2p_{1}^2}+\frac{p^s}{N}\sum_{j=1}^{N-1}|\sum_{i=0}^{p_1^{m-1}-1}\xi_N^{-ji}|^2\nonumber\\
&=& \frac{N(p_{1}-1)^2\ell^2}{w^2p_{1}^2}+\frac{p^s}{N}\sum_{i,k=0}^{p_1^{m-1}-1}\left(\sum_{j=0}^{N-1}\xi_N^{j(k-i)} -1\right)\nonumber\\
&=& \frac{N}{p_1^2}\left(\frac{(p_1-1)^2\ell^2}{w^2}+p^s(p_1-1)\right)\nonumber\\
&=& \frac{N}{p_1}\cdot \frac{(p_1-1)\ell}{w},
\end{eqnarray}
where in the last step of the above computations we used the condition (\ref{cond2}). Let $\kappa=\frac{N}{p_1}\cdot \frac{\ell (p_1-1)}{w}$ and $S=\{a\in \Z/N\Z\mid \varphi(a)\neq 0\}$ (that is, $S$ is the support of $\varphi$). We have
$$0\leq \sum_{a\in S}(\varphi(a)-\epsilon)^2=|S|-\kappa.$$
On the other hand, from $\sum_{a\in G}\varphi(a)^{2}=\kappa$ and $\varphi(a)\in \Z$, we have $\kappa\geq |S|$. Therefore we must have $|S|=\kappa$ and $\sum_{a\in S}(\varphi(a)-\epsilon)^2=0$. Hence $\varphi(a)\in \{0,\epsilon\}$ for all $a\in G$. It follows that $\psi(\gamma^a D)$, $0\leq a\leq N-1$, take only two values. By Theorem~\ref{srgequiv}, Cay$(\F_q,D)$ is an srg. The proof is now complete. $\hfill\Box$

\hspace{0.1in}

\begin{Rem}\label{remark} 
Condition~(\ref{cond1}) is equivalent to 
\begin{eqnarray}\label{cond3}
\frac {p^b (1-p_1)\ell}{w} \equiv \pm 1 \pmod{p_1}
\end{eqnarray}
This can be seen as follows. If we square both sides of (\ref{cond1}), we obtain $\frac{p^{2b+f-\tilde{f}}(1-p_1)^2\ell ^2}{w^2}\equiv 1\pmod {p_1}$. Noting that $p^f\equiv 1\pmod {p_1}$ and $p^{\tilde{f}}\equiv 1\pmod {p_1}$, we have $\frac{p^{2b}(1-p_1)^2\ell ^2}{w^2}\equiv 1\pmod {p_1}$. Since $p_1$ is prime, we must have $\frac {p^b (1-p_1)\ell}{w} \equiv \pm 1 \pmod{p_1}$.  The converse can be proved similarly. We comment that (\ref{cond3}) is much easier to use since it does not involve $m$ any more. 
\end{Rem}

\begin{coro}\label{uppdown}
Let $p_1$ be a prime, $m\geq 1$, $N=p_1^m$. Let $p\neq p_1$ be a prime, $f={\rm ord}_N(p)$, $w=\phi(N)/f$, and $q=p^f$. Assume that $-1\not\in\langle p\rangle$, $\gcd(p(p-1), N)=1$ and $w|(p_1-1)$. Let $\tilde{f}={\rm ord}_{p_1}(p)$, $D$ be defined as in (\ref{DefD}), and $\tilde{D}$ the subgroup of $\F_{p^{\tilde{f}}}^*$ of index $p_1$. Then Cay($\F_{q}$, $D$) is an srg if and only if Cay($\F_{p^{\tilde{f}}}$, $\tilde{D}$) is an srg.
\end{coro}

The proof is clear by Theorem~\ref{srgcond} and the above remark. We omit the details.

\section{New Infinite families of strongly regular Cayley graphs}

In this section, we give three infinite families of srg which are obtained by using Theorem~\ref{eigennumber} and Theorem~\ref{srgcond}.

\noindent{\bf Example 4.1.} Let $p=11$, $p_{1}=43$ and $N=p_{1}^{m}$ for $m\ge 1$. It is easy to use induction to prove that $\order_{43^{m}}(11)=\phi(43^{m})/6$ for all $m\geq 1$. Let $\F_{q}$ be the finite field of order $q=11^{f}$, where $f=\phi(43^{m})/w$, $w=6$. We claim that Cay($\F_{q}$, $D$), with $D=\cup_{i=0}^{p_{1}^{m-1}-1}C_{i}$, is an srg. We could use Corollary~\ref{uppdown} together with the result in Table I to prove this claim. But we prefer to do it without relying on the result in Table I.

In this example, we have $\tilde{f}=7$ and $b=3$ (here $b$ is obtained by computing ${\rm min}\{b_0, b_1, \ldots ,b_{5}\}$ and $b_j=\frac{1}{p_1} \sum _{z \in ([1,p_1-1]\cap \tilde{C}_j)} z$, $0\leq j\leq 5$). It follows that $s=1$. Since $\frac{3}{6}(43-(43-1)\times \frac{3}{6})=11$, (\ref{cond2}) is satisfied with $\ell=3$. Next  $\frac{3\times (1-43)\times 11^{3}}{6}\equiv -1\;(\modulo 43)$, we see that (\ref{cond3}) is satisfied. It follows by Theorem~\ref{srgcond} and Remark~\ref{remark} that Cay($\F_{q}$, $D$) is a strongly regular graph. This family of srg generalizes Example 5 in Table I.

\vspace{0.1in}

\noindent{\bf Example 4.2.} Let $p=5$, $p_{1}=31$ and $N=p_{1}^{m}$ for $m\ge 1$. It is easy to use induction to prove  that $\order_{31^{m}}(11)=\phi(31^{m})/10$. Let $\F_{q}$ be the finite field of order $q=5^{f}$, where $f=\phi(31^{m})/w$ with $w=10$. Now $\tilde{f}=3$. Let $\tilde{D}$ be the subgroup of $\F_{5^3}^*$ of index $p_1=31$. Then $\tilde{D}$ is nothing but $\F_5^*$, i.e. the multiplicative group of the prime subfield of $\F_{5^3}$. Trivially Cay($\F_{p^{\tilde{f}}}$, $\tilde{D}$) is an srg. By Corollary~\ref{uppdown}, Cay($\F_{q}$, $D$) with $D=\cup_{i=0}^{p_{1}^{m-1}-1}C_{i}$ is an srg.

\vspace{0.1in}

\noindent{\bf Example 4.3.} Let $p=2$, $p_{1}=127$ and $N=p_{1}^{m}$ for $m\ge 1$. Again it is easy to use induction to prove that $\order_{127^{m}}(11)=\phi(127^{m})/18$. Let $\F_{q}$ be the finite field of order $q=5^{f}$, where $f=\phi(127^{m})/w$ with $w=18$. Now $\tilde{f}=7$. Let $\tilde{D}$ be the subgroup of $\F_{2^7}^*$ of index $p_1=127$. Then $\tilde{D}$ is nothing but $\F_2^*=\{1\}$. Trivially Cay($\F_{p^{\tilde{f}}}$, $\tilde{D}$) is an srg. By Corollary~\ref{uppdown}, Cay($\F_{q}$, $D$) with $D=\cup_{i=0}^{p_{1}^{m-1}-1}C_{i}$ is an srg. It should be noted that as an srg,  Cay($\F_{q}$, $D$) is not trivial at all.

\section{Concluding remarks}

In conclusion, we generalize the construction in \cite{GXY} to give new infinite families of srg. The main new result here is Theorem~\ref{srgcond}, which gives necessary and sufficient conditions for $D$ in (\ref{DefD}) to give rise to srg.

After finishing the research of this paper, we became aware of the very interesting paper \cite{momihara2012} by Momihara. In \cite{momihara2012}, 
the author gave a recursive construction of strongly regular Cayley graphs, generalizing all 11 sporadic examples in the statement of the 
Schmidt-White conjecture into infinite families. In particular, the two index 6 examples are generalized into infinite families while 
we could only generalize one of the index 6 examples in this paper. However, the approach taken in our paper is different from that of \cite{momihara2012} 
since ours is a direct construction. Also we obtained two conditions (\ref{cond1}) and (\ref{cond2}) which are necessary and sufficient for our construction 
to give rise to an srg. These conditions reveal an interesting connection between strongly regular Cayley graphs and cyclic difference sets 
in $(\Z/p_1\Z, +)$, which will be useful in future investigation of strongly regular Cayley graphs and cyclic difference sets.


\begin{thebibliography}{99}
\bibitem{AL} Aubry, Y., Langevin, P., On the weight of binary irreducible cyclic codes, Workshop on Coding and Cryptography WCC'05 ( Springer )  {\bf 3969}  Norway, 46--54 ( 2006 ).

\bibitem{B} Bose, R. C.: Strongly regular graphs, partial geometries, and partially balanced designs, Pacific J. Math. {\bf 13} (1963), 389–419.


\bibitem{BEW} Berndt, B. C., Evans, R. J., Williams, K. S.: {\it Gauss and Jacobi Sums}, A Wiley-Interscience Publication, New York, (1998).


\bibitem{BH} Brouwer,A. E., Haemers, W. H.: {\it Spectra of Graphs}, Universitext, Springer-Verlag, New York, (2012).

\bibitem{BMW} Baumert, L.D., Mills, W. H., Ward, R. L.: Uniform Cyclotomy, Journal of Number Theory {\bf 14}, 67-82 (1982).

%\bibitem{BWX} Brouwer, A. E., Wilson, R. M., Xiang, Q.: Cyclotomy and Strongly Regular Graphs, J. Alg. %Combin. {\bf 10}, 25--28 (1999).

\bibitem{CK} Calderbank, R., Kantor, W.M.: The geometry of two-weight codes, Bull. Lond. Math. Soc. {\bf 18 (2)}, 97-22 (1986).

\bibitem{DL} De Lange, C.L.M.: Some new cyclotomic strongly regular graphs, J. Algebraic Combin. {\bf 4}, 329-330 (1995).

\bibitem{FYL} Feng, K., Yang, J., and Luo, S.: Gauss sum of index 4: (1) Cyclic case, Acta Math. Sin. (English Ser.) {\bf 21-6}, 1425-1434 (2005).

\bibitem{FX} Feng, T., Xiang, Q.: Strongly regular graphs from unions of cyclotomic classes, J. Combin. Theory (B), in press, http://dx.doi.org/10.1016/j.jctb.2011.10.006.

\bibitem{FMX} Feng, T., Momihara, K., Xiang, Q., Constructions of strongly regular Cayley graphs and skew Hadamard difference sets from cyclotomic classes, available at {\tt arXiv:1201.0701}.

\bibitem{GXY} Ge, G., Xiang, Q. and Yuan, T.: Constructions of Strongly regular Cayley graphs using index 4 Gauss sums, J. Alg. Combin., in press, DOI 10.1007/s10801-012-0368-y.

\bibitem{GR} Godsil, C., Royle, G.: Algebraic Graph Theory, Springer-Verlag, New York, Inc, (2001).

\bibitem{IM} Ikuta, T., Munemasa, A.: Pseudocyclic association schemes and strongly regular graphs, European J. Combin. {\bf 31}, 1513-1519 (2010).

\bibitem{IR} Ireland, K., Rosen, M.: A classical introduction to modern number theory 2nd Ed., Springer-Verlag, New York, (1990).

%\bibitem{Lang} Lang, S.: Cyclotomic Fields, I and II, 2nd Ed., Springer-Verlag, New York, (1990).

\bibitem{L} Langevin,P.: A new class of two-weight codes, in: S. Cohen, H. Niederreiter (Eds.), Finite Fields and Applications, Glasgow, 1995, in: London Math. Soc. Lecture Note Ser., vol. {\bf 233}, Cambridge University Press, 1996, pp. 181-187.

\bibitem{Ma} Ma, S. L., A survey of partial difference sets, Des. Codes Cryptogr. {\bf 4}, 221-261 (1994).

\bibitem{momihara2012} Momihara, K.: Strongly regular Cayley graphs, skew Hadamard difference sets, and rationality of relative Gauss sums, available at {\tt arXiv:1202.6414}. 

\bibitem{T} Thaine, F.: Properties that characterize Gaussian periods and cyclotomic numbers, Proc. AMS {\bf 124},  35-45 (1996).

\bibitem{SW} Schmidt, B., White, C.: All two-weight irreducible cyclic codes?, Finite Fields and Their Appl. {\bf 8}, 1-17 (2002).

\bibitem{VW} Van Lint, J. H., Wilson, R.: Combinatorics, 2nd Ed, Cambridge Press, (2001).

\end{thebibliography}
\end{document}